\documentclass{amsart}

\newtheorem{theorem}{Theorem}[section]
\newtheorem{lemma}[theorem]{Lemma}
\newtheorem{corollary}[theorem]{Corollary}
\newtheorem{proposition}[theorem]{Proposition}

\theoremstyle{definition}

\newtheorem{example}[theorem]{Example}

\theoremstyle{remark}

\numberwithin{equation}{section}

\newcommand{\abs}[1]{\lvert#1\rvert}

\newcommand{\im}{\mathop{\mathrm{Im}}\nolimits}
\newcommand{\re}{\mathop{\mathrm{Re}}\nolimits}

\newcommand{\bbC}{{\mathord{\mathbb C}}}
\newcommand{\bbF}{{\mathord{\mathbb F}}}
\newcommand{\bbQ}{{\mathord{\mathbb Q}}}
\newcommand{\bbR}{{\mathord{\mathbb R}}}
\newcommand{\bbT}{{\mathord{\mathbb T}}}
\newcommand{\bbD}{{\mathord{\mathbb D}}}
\newcommand{\bbZ}{{\mathord{\mathbb Z}}}

\newcommand{\al}{{\mathord{\alpha}}}
\newcommand{\be}{{\mathord{\beta}}}
\newcommand{\De}{{\mathord{\Delta}}}

\newcommand{\Ga}{{\mathord{\Gamma}}}

\newcommand{\si}{{\mathord{\sigma}}}
\newcommand{\ka}{{\mathord{\kappa}}}

\newcommand{\la}{{\mathord{\lambda}}}

\newcommand{\thh}{{\mathord{\theta}}}

\newcommand{\eps}{{\mathord{\epsilon}}}
\newcommand{\ze}{{\mathord{\zeta}}}
\newcommand{\om}{{\mathord{\omega}}}

\newcommand{\cM}{\mathord{\mathcal{M}}}
\newcommand{\cP}{\mathord{\mathcal{P}}}

\newcommand{\fra}{{\mathord{\mathfrak{a}}}}
\newcommand{\frb}{{\mathord{\mathfrak{b}}}}
\newcommand{\frc}{{\mathord{\mathfrak{c}}}}
\newcommand{\frd}{{\mathord{\mathfrak{d}}}}

\newcommand{\frg}{{\mathord{\mathfrak{g}}}}
\newcommand{\frl}{\mathord{\mathfrak{l}}}
\newcommand{\frk}{{\mathord{\mathfrak{k}}}}
\newcommand{\frm}{{\mathord{\mathfrak{m}}}}
\newcommand{\frs}{{\mathord{\mathfrak{s}}}}

\newcommand{\frt}{{\mathord{\mathfrak{t}}}}
\newcommand{\frv}{{\mathord{\mathfrak{v}}}}
\newcommand{\frz}{{\mathord{\mathfrak{z}}}}
\newcommand{\frN}{{\mathord{\mathfrak{N}}}}

\newcommand{\sfh}{{\mathord{\mathsf{h}}}}
\newcommand{\sfe}{{\mathord{\mathsf{e}}}}
\newcommand{\sff}{{\mathord{\mathsf{f}}}}
\newcommand{\sfk}{{\mathord{\mathsf{k}}}}

\newcommand{\GL}{\mathop{\mbox{\rm GL}}\nolimits}
\newcommand{\BL}{\mathop{\mbox{\rm L}}\nolimits}

\newcommand{\OO}{\mathop{\mbox{\rm O}}\nolimits}

\newcommand{\SU}{\mathop{\mathrm{SU}}\nolimits}

\newcommand{\rU}{\mathop{\mathrm{U}}\nolimits}
\newcommand{\ad}{\mathop{\mathrm{ad}}\nolimits}
\newcommand{\Ad}{\mathop{\mathrm{Ad}}\nolimits}
\newcommand{\Aut}{\mathop{\mathrm{Aut}}\nolimits}
\newcommand{\conv}{\mathop{\mathrm{conv}}\nolimits}
\newcommand{\cone}{\mathop{\mathrm{cone}}\nolimits}
\newcommand{\codim}{\mathop{\mathrm{codim}}\nolimits}

\newcommand{\spec}{\mathop{\mathrm{spec}}\nolimits}
\newcommand{\spann}{\mathop{\mathrm{span}}\nolimits}

\newcommand{\sll}{\mathop{\mathrm{sl}}\nolimits}
\newcommand{\su}{\mathop{\mathrm{su}}\nolimits}
\newcommand{\scal}[2]{\left<#1,#2\right>}
\newcommand{\clos}{\mathop{\mbox{\rm clos}}}
\newcommand{\card}{\mathop{\mbox{\rm card}}}

\newcommand{\csn}{{(\bbC^*)^n}}
\newcommand{\one}{{\mathord{\bf 1}}}

\newcounter{xxx}

\let\wh=\widehat
\let\wt=\widetilde

\begin{document}

\title{Orbits of tori extended by finite groups and their
polynomial hulls:  the case of connected complex orbits}

\author{V.M.~Gichev}
\address{Omsk Branch of Sobolev Institute of Mathematics,
Pevtsova 13, 644099, Omsk, Russia}
\email{gichev@ofim.oscsbras.ru}
\thanks{The author was partially supported by RFBR Grants 06-08-01403 and
06-07-89051.}
\subjclass{Primary 32E20; 
Secondary 32M15, 32M05. }
\date{}

\keywords{Polynomial hulls, bounded symmetric domains}

\begin{abstract}
Let $V$ be a complex linear space, $G\subset\GL(V)$ be a compact
group. We consider the problem of description of polynomial hulls
$\wh{Gv}$ for orbits $Gv$, $v\in V$, assuming that the identity
component of $G$ is a torus $T$. The paper contains a universal
construction for orbits which satisfy the inclusion $Gv\subset
T^\bbC v$ and a characterization of pairs $(G,V)$ such that it is
true for a generic $v\in V$. The hull of a finite union of
$T$-orbits in $T^\bbC v$ can be distinguished  in $\clos T^\bbC v$
by a finite collection of inequalities of the type
$\abs{z_1}^{s_1}\dots\abs{z_n}^{s_n}\leq c$. In particular, this
is true for $Gv$. If powers in the monomials are independent of
$v$, $Gv\subset T^\bbC v$ for a generic $v$, and either the center
of $G$ is finite or $T^\bbC$ has  an open orbit, then the space
$V$ and the group $G$ are products of standard ones; the latter
means that $G=S_nT$, where $S_n$ is the group of all permutations
of coordinates and $T$ is either $\bbT^n$ or $\SU(n)\cap\bbT^n$,
where $\bbT^n$ is the torus of all diagonal matrices in $\rU(n)$.
The paper also contains a description of polynomial hulls for
orbits of isotropy groups of bounded symmetric domains. This
result is already known, but we formulate it in a different form
and supply with a shorter proof.
\end{abstract}

\maketitle 
\section*{Introduction}
Let $V$ be a finite-dimensional complex linear space and
$G\subset\GL(V)$ be a compact subgroup of $\GL(V)$. We consider
the problem of description of polynomially convex hulls for orbits
$O_v=Gv$, $v\in V$. The {\it polynomially convex hull} (or {\it
polynomial hull}) $\wh Q$ of a compact set $Q\subset V$ is defined
as
\begin{eqnarray}\label{defph}
\wh Q=\{z\in V:\,\abs{p(z)}\leq\sup_{\ze\in
Q}|p(\ze)|\quad\mbox{for all}\quad p\in\cP(V)\},
\end{eqnarray}
where $\cP(V)$ is the algebra of all holomorphic polynomials on
$V$. It is usually difficult to find $\wh Q$. For $Q=Gv$, the
answer is known if $G$ is an isotropy group of a bounded symmetric
domain in $\bbC^n$. Paper \cite{Ka} contains a description of
$G$-invariant polynomially convex compact sets, including hulls of
orbits ($Q\subset V$ is {\it polynomially convex} if $\wh Q=Q$);
it continues paper \cite{KZ} and uses results of \cite{FB}. On the
other hand, it is known that an orbit of a compact linear group is
polynomially convex if and only if the complex orbit $G^\bbC v$ is
closed and $Gv$ is its real form (\cite{GL}). The cases
$G=\rU(2),\SU(2)$ were considered in \cite{An}, \cite{DG}. The
problem of determination of polynomial hulls of orbits admits the
following natural generalization: given a homogeneous space $M$ of
a compact group $G$, describe maximal ideal spaces $\cM_A$ of
$G$-invariant closed subalgebras $A$ of $C(M)$, where $C(M)$ is
the Banach algebra of all continuous complex-valued functions on
$M$ with the sup-norm. If $A$ is generated by a finite-dimensional
invariant subspace, then $\cM_A$ can be realized as the polynomial
hull of an orbit. Paper \cite{Gi} contains a description of
$\cM_A$ for bi-invariant algebras on compact groups and partial
results on spherical homogeneous spaces. Maximal ideal spaces for
$\rU(n)$-invariant algebras on spheres in $\bbC^n$ are described
in \cite{Kan}.

In this paper we consider orbits $Gv$ of groups $G=FT$, where
$F\subseteq G$ is a finite subgroup and $T$ is a torus, such that
$G^\bbC v=T^\bbC v$. Let $\frt\subseteq\frg\frl(V)$ be the Lie
algebra of $T$ and set $\frt^\bbR=i\frt$,
$T^\bbR=\exp(\frt^\bbR)$. Suppose that $v\in V$ has a trivial
stable subgroup in $T$ and let $X\subset T^\bbR v$ be finite. The
hull of $Y=TX$ admits a simple description. If $X=\{v\}$, then
$\wh{Y}=\wh{Tv}$ is the closure of $T\exp(C_T)v$, where $C_T$ is a
cone in $\frt^\bbR$. If $T^\bbC$ is closed, then $\wh
Y=T\exp(Q_X)v$, where $Q_X\subseteq\frt^\bbR$ is a convex polytope
(the convex hull of the inverse image of $X$ for the mapping
$\xi\to\exp(\xi)v$, $\xi\in\frt^\bbR$). Any segment in $Q_X$
corresponds to an analytic strip or an annulus in $\wh Y$. In
general, $\wh Y$ is the union of $\wh {Tu}$, where $u$ runs over
$\exp(Q_X)v$. Also, $\wh Y$ is distinguished in $\clos T^\bbC v$
by a finite family of monomial inequalities of the type
\begin{eqnarray}\label{moine}
\abs{z_1}^{s_1}\dots\abs{z_n}^{s_n}\leq c,
\end{eqnarray}
where $c\geq0$ and $s=(s_1,\dots,s_n)\in\bbR^n$ depend on $v$ and
$X$. Vectors $s$ correspond to normals of faces of $C_T+Q_X$.

Thus, the problem of determination of $\wh{Gv}$ is not difficult
if $Gv\subset T^\bbC v$. The latter is equivalent to the
assumption that the complex orbit $G^\bbC v$ is connected. In
Example~\ref{maico}, we give a construction for orbits which
satisfy this condition; here is a sketch. The group $G=FT$ acts on
the space $V=C(K)$, where $K$ is a finite $F$-invariant subset of
$\frt^*$: $F$ acts naturally on $C(K)$, $\frt=\frt^{**}$ is
naturally embedded into $C(K)$, and $T=\exp(\frt)$ acts on $C(K)$
by multiplication. If $v\in C(K)$ is an $F$-invariant function,
then $Gv\subset T^\bbC v$. According to Theorem~\ref{isase}, each
connected complex orbit can be realized in this way. Further, we
describe pairs $(V,G)$ such that
\begin{eqnarray}\label{genev}
Gv\subset T^\bbC v\quad\mbox{for a generic}\quad v\in V.
\end{eqnarray}
By Theorem~\ref{charf}, under the additional assumption that the
complex linear span of $T^\bbC v$ coincides with $V$, this happens
if and only if the group $G^\bbC Z$, where $Z$ is the centralizer
of $G$ in $\GL(V)$, has an open orbit in $V$. There are two
extreme cases: (A) $Z\subseteq G^\bbC$; (B) $G$ has a finite
center. An example for (A) is the group $G=S_n\bbT^n$ acting in
$\bbC^n$, where $\bbT^n$ is the torus of all diagonal matrices in
$\rU(n)$ and $S_n$ is the group of all permutations of
coordinates. Replacing $\bbT^n$ with $\SU(n)\cap\bbT^n$, we get an
example for (B). Example~\ref{maict} contains a construction for
pairs $(V,G)$ that satisfy (\ref{genev}). Theorem~\ref{asexa}
states that the construction is a universal one. In
Theorem~\ref{huini}, we determine pairs which satisfy
(\ref{genev}) and the following condition:
\begin{eqnarray}\label{coxet}
\mbox{vectors}\ s\ \mbox{in}\ (\ref{moine})\ \mbox{are independent
of}\ v.
\end{eqnarray}
 The paper also contains a description of
hulls $\wh{Gv}$ for $G=\Aut_0(D)$, where $D$ is a bounded
symmetric domain in the canonical realization and $\Aut_0(D)$ is
the stable subgroup of zero, which coincides with the group of all
linear automorphisms of $D$. These hulls have already been
described: the final step was done in paper \cite{Ka}, which
essentially used \cite{KZ}, partial results appear in \cite{Sa}
and \cite{FB}. Most of them use the technique of Jordan triples
and Jordan algebras. We use Lie theory, in particular, an explicit
construction of paper \cite{Wo} for a maximal abelian subspace
$\fra$. A compact group acting in a Euclidean space is called {\it
polar} if there exists a subspace (a {\it Cartan subspace}) such
that each orbit meets it orthogonally. The group $G$ is  polar in
the ambient linear space $\frd$, and $\fra$ is the Cartan subspace
for $G$. Real polar representations are classified in paper
\cite{Da}; they are orbit equivalent (i.e., have the same orbits)
to isotropy representations of Riemannian symmetric spaces. If $D$
is a polydisc $\bbD^n\subset\bbC^n$, where $\bbD$ is the unit disc
in $\bbC$, then $G=S_n\bbT^n$; the polynomial hulls $\wh{Gv}$ are
determined by the inequalities
\begin{eqnarray}\label{mukin}
\mu_k(z)\leq\mu_k(v),
\end{eqnarray}
where $k=1,\dots,n$ and $\mu_k$ are defined by
\begin{eqnarray}\label{mukde}
\mu_k(z)=\max\{\abs{{z_{\si(1)}}\dots{z_{\si(k)}}}:\,\si\in S_n\}.
\end{eqnarray}
The general case can be reduced to this one in the following way.
Any bounded symmetric domain $D\subset\frd$ of rank $n$ admits an
equivariant embedding of $\bbC^n$ to $\frd$, which induces an
embedding of $\bbD^n$ to $D$, such that $\bbR^n\subset\bbC^n$ is
the maximal abelian subspace $\fra$, and, for any  $v\in\fra\,$,
the hull of $\Aut_0(D)v$ is the orbit of the hull of
$\Aut_0(\bbD^n)v$. Each $\mu_k(z)$ has a unique continuation to a
$K$-invariant function on $\frd$. The extended functions determine
hulls by the same inequalities. Moreover, they are
plurisubharmonic and can be treated as products of singular values
of $z\in\frd$ or as norms of exterior powers of adjoint operators
in suitable spaces. The subsystem of long roots of the restricted
root system (i.e., the root system for $\fra$) has type $nA_1$;
this defines the above embedding $\bbC^n\to\frd$. Furthermore,
this makes it possible to determine hulls in terms of the adjoint
representation (Theorem~\ref{last}). Thus, there is no need to
consider different types of domains separately.

The reduction to the case of a torus extended by a finite group,
which is described above, is contained in Section~\ref{isohe} (in
papers \cite{Ka}, \cite{Sa}, the problem is also reduced to this
case by another method). It does not use essentially the results
of the previous sections (only Proposition~\ref{autdn}, in proof
of Theorem~\ref{last}). These extensions satisfy conditions
(\ref{genev}) and (\ref{coxet}); in addition, they possess the
property that the complexified groups have open orbits. According
to Theorem~\ref{huini}, any group with these properties is the
product of groups $S_n\bbT^n$ acting in $\bbC^n$; it admits a
natural realization as a group of automorphisms of a bounded
symmetric domain (Corollary~\ref{conop}).

The following simple examples illustrate the case $Gv\not\subseteq
T^\bbC v$ and  show that condition~(\ref{genev}) is essential. Let
$G=S_n\bbT^n$, and let $\eps_1,\dots,\eps_n$ be the standard base
in $\bbC^n$. Then $\wh{G\eps_1}$ is the closure the union of discs
$\bbD\eps_k$, $k=1,\dots,n$. Set $H=S_n\bbT$, where $\bbT$ acts by
$z\to e^{it}z$, $t\in\bbR$, $z\in\bbC^n$. Then
$\wh{H\eps_1}=\wh{G\eps_1}$. For $v=\eps_1+\eps_2$, $\wh{Gv}$ is
the closure of the union of $n\choose2$ bidiscs but $\bbT^n$
contains no proper torus $T$ such that $\wh{Gv}=\wh{Hv}$ for
$H=S_nT$. However, for any subgroup $F\subseteq S_n$ which acts
transitively on 2-sets and $H=F\bbT^n$ we have $\wh{Gv}=\wh{Hv}$.

\section{Preliminaries}
We keep the notation of Introduction, in particular, (\ref{defph})
and (\ref{mukde}). Linear spaces are supposed to be finite
dimensional and complex unless the contrary is explicitly stated.
"Generic" means "in some open dense subset". Throughout the paper,
we use the following notation:
\begin{itemize}
\item[] $\bbD$ and $\bbT$ are the open unit disc and the unit
circle in $\bbC$, respectively; \item[] $V$ denotes a complex
linear space (except for Section~\ref{isohe}); \item[] if $V$ is
equipped with a linear base identifying it with $\bbC^n$, then
$\bbT^n$ is the group of all diagonal  unitary transformations;
\item[] $\bbZ_2^n$ consists of all transformations in $\bbT^n$
with eigenvalues $\pm1$;  \item[] $\eps_1,\dots,\eps_n$ is the
standard base in $\bbC^n$ and $\bbR^n$;  \item[] $\bbR^n_+$ is the
set of vectors in $\bbR^n$ with positive entries; \item[] $S_K$
denotes the group of all permutations of a finite set $K$; if
$K=\{1,\dots,n\}$, then $S_K=S_n$; \item[] $C(K)$ is the algebra
of all complex-valued functions on $K$; \item[] $\one$ is the
identity of $C(K)$; \item[] $G\subset\GL(V)$ is a compact group
whose identity component is a torus $T$ (except for
Section~\ref{isohe}); \item[] $\frt\subset\frg\frl(V)$ is the Lie
algebra of $T$, $\frt^\bbR=i\frt$, $\frt^\bbC=\frt+\frt^\bbR$;
\item[] $T^\bbR=\exp(\frt^\bbR)$, $T^\bbC=\exp(\frt^\bbC)$;
\item[] $\bbC^*=\bbT^\bbC=\bbC\setminus \{0\}$; \item[] $\check
T=\mathop{\mbox{\rm{Hom}}}(T,\bbT)$ is the dual group to $T$;
\item[] $\Aut(D)$ is the group of all holomorphic automorphisms of
a domain $D\subset V$, $\Aut_0(D)=\Aut(D)\cap\GL(V)$;  \item[]
$\cone X$ denotes the least convex cone which contains the set
$X$;\item[] $\conv X$ is the convex hull of $X$; \item[] $\clos X$
is the closure of $X$; \item[] $\spann_\bbF X$ is the linear span
of $X$ over the field $\bbF=\bbC,\bbR,\bbQ$.

\end{itemize}
Clearly, $\exp$ is bijective on $\frt^\bbR$ and $T^\bbR\cong
T^\bbC/T$. The differentiating at the identity $e$ defines an
embedding of $\check T$ into the dual space $\frt^*$: $\chi\to-i
d_e\chi$, where $\chi\in\check T$. This is a lattice in the vector
group $\frt^*$, moreover, $T\cong\frt/L$, where $L$ is the dual
lattice to $\check T$ in $\frt$. For $\chi\in\check T$, let
$$V_\chi=\{v\in V:\,gv=\chi(g)v\ \mbox{\rm for all}\ g\in T\}$$
be the corresponding isotypical component of $V$. Then
\begin{eqnarray}\label{decv}
V=\sum\nolimits_{\chi\in\check T}\oplus\, V_\chi.
\end{eqnarray}
We assume that $V$ is equipped with a $G$-invariant inner product
$\scal{\ }{\ }$. Then decomposition (\ref{decv}) is orthogonal.
Let $\spec(v)$ denote the spectrum of $v\in V$ (the set of
$\chi\in\check T$ such that the $\chi$-component of $v$ is
nonzero); for $X\subseteq V$,
$$\spec(X)=\cup_{x\in X}\spec(x).$$
We say that $T$ has a {\it simple spectrum} if
\begin{eqnarray}\label{sispe}
\dim V_\chi\leq1
\end{eqnarray}
for all $\chi\in\check T$. If (\ref{sispe}) is true, then there
exists a unique (up to scaling factors) orthogonal base in $V$
which agree with (\ref{decv}) and a unique maximal torus $\bbT^n$
in $\GL(V)$ which contains $T$. In what follows, we assume that
(\ref{sispe}) holds; we shall see in the next section that such
assumption is not restrictive. Thus, we may fix an identification
\begin{eqnarray}\label{vcnck}
V=\bbC^n=C(K),
\end{eqnarray}
where $K=\{1,\dots,n\}$. If $F$ is a subgroup of $S_K$, then
$C(K)^F$ denotes the set of all $F$-invariant functions on $K$;
clearly, $\one\in C(K)^F$. Further, $\csn$ is the multiplicative
group of all invertible functions in $C(K)$, $\bbT^n$ consists of
functions with values in $\bbT$, and $(\bbT^n)^\bbC=\csn$. The Lie
algebra of $\bbT^n$  is realized as $i\bbR^n\subset\bbC^n$. The
embedding $T\to\bbT^n$ induces embeddings of the Lie algebra and
the fundamental group: $\frt\to i\bbR^n$, $\pi_1(T)\to
i\bbZ^n\subset i\bbR^n$, respectively. Let $\Ga$ be the image of
$\pi_1(T)$. Then $\spann_\bbR\Ga=\frt$; moreover, $\frt\cap
i\bbZ^n=\Ga$ and $\frt/\Ga=T$. The dual mapping
$\check\bbT^n\to\check T$, which is defined by the restriction of
characters $e^{-i\scal{x}{y}}$, where $x\in i\bbZ^n$,  to $\frt$,
is the orthogonal projection $\pi_\frt:\,i\bbZ^n\to\frt$. Thus,
$\Ga$ is a subgroup of finite index in $\check T=\pi_\frt
i\bbZ^n$. Vectors in $\spann_\bbQ\check T$ are called {\it
rational}. The image of $\frt$ in $i\bbR^n$ can be distinguished
by linear equations with integer coefficients. Hence,
$\clos(T^\bbC v)$, for a generic $v\in V$, is the set of all
solutions to a finite number of equalities with holomorphic
monomials. Thus,
$Y\subset T^\bbC$ implies
$\wh{Y}\subset\clos(T^\bbC v)$.
Set
\begin{eqnarray}\label{defct}
C_T=\frt^\bbR\cap\clos(-\bbR^n_+),
\end{eqnarray}
The cone $iC_T$ is dual to $\cone(\spec V)\subseteq\frt^*\subseteq
i\bbR^n$. If $-\xi\in\clos(\bbR^n_+)$, then
$\iota=\lim_{t\to+\infty}\exp(t\xi)$ is an idempotent in $C(K)$
such that the multiplication by the complementary idempotent
$\one-\iota$ is a projection onto $\spann_\bbC(\spec(\xi))$. Set
\begin{eqnarray}\label{defit}
I_T=\{\lim\nolimits_{t\to+\infty}\exp(t\xi):\,\xi\in C_T\}.
\end{eqnarray}
Clearly, $I_T$ is finite and contains $\one$.
\begin{lemma}\label{cloct}
The closure of $\exp(C_T)$ is equal to $I_T\exp(C_T)$.
\end{lemma}
\begin{proof}
Due to the evident inclusion $\clos(\exp(C_T))\supseteq
I_T\exp(C_T)$, it is sufficient to prove that the set
$S_T=I_T\exp(C_T)$ is closed. Clearly, $S_T$ is an abelian
semigroup. The cone $C_T$ is polyhedral; hence, it is finitely
generated:
\begin{eqnarray*}
C_T=\cone\{\xi_1,\dots,\xi_m\},
\end{eqnarray*}
where $\bbR^+\xi_k$ are the extreme rays of $C_T$, $k=1,\dots,m$.
Obviously, $I_T$ is a finite semigroup, which is generated by the
idempotents $\lim_{t\to+\infty}\exp(t\xi_k)$. Thus, the
correspondence
$(e^{-t_1},\dots,e^{-t_m})\to\exp(t_1\xi_1+\dots,t_m\xi_m)$
defines a mapping of $(0,1]^m$ onto $\exp(C_T)$, which
continuously extends to $[0,1]^m$. It follows that its image is
closed and coincides with $S_T$.
\end{proof}
Note that there is a natural one-to-one correspondence between
$I_T$ and the set of faces of $C_T$.

\section{Hulls of finite unions of $T$-orbits in a $T^\bbC$-orbit}
Let $v\in T^\bbC$. If $v=\sum_{\chi\in\check T}v_\chi$, where
$v_\chi\in V_\chi$,  $g\in T^\bbC$, and $u=gv$, then
$u=\sum_{\chi\in\check T}\chi(g)v_\chi$. Since $\chi(g)\neq0$ for
all $g\in G$ and $\chi\in\check T$, we get
\begin{eqnarray}
&u\in T^\bbC v\quad\Longrightarrow\quad\spec(u)=\spec(v);\label{speceq}\\
&\dim\left(V_\chi\cap\spann_\bbC Tv\right)\leq1\quad\mbox{\rm for
all}\quad v\in V\quad\mbox{\rm and}\quad\chi\in\check
T.\label{simsp}
\end{eqnarray}
Thus, the assumption that $T$ has a simple spectrum in $V$ is not
restrictive in the problem of description of polynomial hulls of
orbits $Gv$ such that $Gv\subset T^\bbC v$. Clearly,
$\bbD^n=\wh{\bbT^n}$ in $\BL(V)$. For each $x\in C_T$ and any
polynomial $p$ on $\BL(V)$, the holomorphic function
$f(\ze)=p(\exp(\ze x))$ is bounded in the halfplane $\Pi:\,\re
\ze\geq0$. Hence, $\exp(\Pi)$ is contained in $\wh T$. On the
other hand, if $z\in \bbD^n\cap T^\bbC$, then $z=t\exp(x)$ for
some $t\in T$ and $x\in C_T$ (the polar decomposition). By
Lemma~\ref{cloct},
\begin{eqnarray*}
\wh T=\clos(\bbD^n\cap T^\bbC)=T\clos(\exp (C_T))=TI_T\exp(C_T).
\end{eqnarray*}
If $v\in\csn$, then $\csn v=\csn$, and the mapping $z\to zv$ is a
linear nondegenerate transformation of $\bbC^n$. Therefore,
\begin{eqnarray}\label{hltgv}
v\in \csn\quad\Longrightarrow\quad\wh{Tv}=\wh Tv=TI_T\exp (C_T)v.
\end{eqnarray}
For an arbitrary $v\in V=\bbC^n$, set
\begin{eqnarray*}
C_T^v=\{\xi\in\frt^\bbR:\,\xi_k\leq0\ \mbox{if}\ v_k\neq0,\
k=1,\dots,n\}.
\end{eqnarray*}
Applying (\ref{hltgv}) to $\spann_\bbC(\spec(v))=\bbC^n v$, we get
\begin{eqnarray}\label{hultv}
\wh{Tv}=T\clos(\exp(C_T^v)v.
\end{eqnarray}
Clearly, $C_T^v$ depends only on $\spec(v)$. For $s\in\bbR^n$ and
$z\in(\bbC^*)^n$, set
\begin{eqnarray*}
\nu_s(z)=\prod\nolimits_{k=1}^n|z_k|^{s_k}.
\end{eqnarray*}
If $s_k\geq0$, then the $k$-th factor in $(\bbC^*)^n$ can be
replaced with $\bbC$ (i.e., $\nu_s$ continuously extends to this
product).

It is well known that for any holomorphically convex $T$-invariant
set $U\subseteq T^\bbC$, the set $\log(U\cap
T^\bbR)\subseteq\frt^\bbR$ is convex. In particular, this is true
for sets of $g\in T^\bbC$ such that $gv\in\wh{TX}$, where
$X\subset T^\bbC v$, $v\in V$. Nevertheless, it is convenient to
have an explicit construction of an analytic strip (or an annulus,
if it is periodic) in a $T^\bbC$-orbit, which corresponds to a
segment that joins two points in $\frt^\bbR$; it is contained in
the following lemma.  Set
\begin{eqnarray*}
S=\{z\in\bbC:\,0\leq\im z\leq 1\}.
\end{eqnarray*}

\begin{lemma}\label{strih}
Let $v\in\bbC^n$ and $u\in T^\bbR v$. Then, there exists
$\xi\in\frt^\bbR$ such that
\begin{eqnarray}\label{anstdi}
\la(z)=\exp(z\xi)v
\end{eqnarray}
is a holomorphic mapping  $\la:\,S\to T^\bbC v$ which satisfies
conditions
\begin{eqnarray*}
\la(\partial S)\subseteq Tv\cup Tu,\\
\la(0)=v,\quad\la(1)=u.
\end{eqnarray*}
If the stable subgroup of $v$ in $T^\bbR$ is trivial, then $\xi$
is unique.
\end{lemma}
\begin{proof} These properties hold for $\xi\in\frt^\bbR$ such
that $\exp(\xi)v=u$; such a $\xi$ exists, since $\exp$ is a
bijection $\frt^\bbR\to T^\bbR$. The last assertion is clear.
\end{proof}

If $\xi\in C_T$, then (\ref{anstdi}) defines an analytic halfplane
in $\wh{Tv}$; for $\Ga$-rational $\xi$, $\la$ is periodic and
defines an analytic disc in $\wh {Tv}$. Together with
Lemma~\ref{strih} this gives a characterization of hulls for
finite unions of $T$-orbits in $T^\bbC$. Suppose that $X\subset
T^\bbR v$ is finite and the stable subgroup of $v$ in $T$ is
trivial. Then, the inverse to the mapping $x\to\exp(x)v$ is well
defined. Let us denote it by $\log_v$ and set
\begin{eqnarray}
&Q_X=\conv(\log\nolimits_v X),\label{defqx}\\
&P_X=Q_X+C_T.\label{defpx}
\end{eqnarray}
The set $P_X$ is a convex polyhedron, which is unbounded if
$C_T\neq0$. Hence, there exists a finite set $\frN_X\subset\bbR^n$
and, for each $s\in\frN_X$, real numbers $c_s$ such that
\begin{eqnarray}\label{depx}
P_X=\{x\in\frt^\bbR:\,\scal{x}{s}\leq c_s\ \mbox{\rm for all}\
s\in\frN_X\}.
\end{eqnarray}
The set $\frN_X$ consists of vectors orthogonal to faces of $P_X$,
whose projections into $\spann_\bbR P_X$ look outside of it;
clearly, it is not unique in general.
\begin{proposition}\label{hulfi}
Let $v\in \csn$. Suppose that $Y\subset T^\bbC v$ is a finite
union of $T$-orbits (including $Tv$), and set $X=T^\bbR v\cap Y$.
Then $X$ is finite and
\begin{eqnarray}
\wh{Y}&=&\clos\left(T\exp(P_X)v\right)\label{grofu}\\
&=& \clos\left\{z\in\csn:\,\nu_s(z)\leq e^{c_s}\nu_s(v), \
s\in\frN_X\right\}\label{inehu}\\
&=& \bigcup\nolimits_{u\in\exp(Q_X)v}\wh{Tu}\label{cuphu}\\
&=&T\exp(P_X)I_Tv,\label{addfa}
\end{eqnarray}
where $Q_X,P_X,\frN_X$ are as above and $I_T$ is defined in
{\rm(\ref{defit})}.
\end{proposition}
\begin{proof}
Due to the polar decomposition, the set $Tu\cap T^\bbR v$, for
each $u\in T^\bbC v$, is nonvoid and consists of a single point.
Hence, $X$ is finite and $Y=TX$. The inclusion $(\exp
Q_X)v\subseteq\wh{Y}$ follows from Lemma~\ref{strih} and
Phragm\'en--Lindel\"of Principle. The inclusion
$\exp(C_T)u\subseteq\wh{Tu}$ is true for any $u\in \bbC^n$. Since
it holds for all $u\in  T\exp(Q_X)$, the left-hand side of
(\ref{grofu}) includes the right-hand side. If $z=\exp(\xi)v$,
where $\xi\in\frt^\bbC$, then $z_k=e^{\xi_k}v_k$, $k=1,\dots,n$;
due to (\ref{depx}), this implies that the right-hand side of
(\ref{grofu}) coincides with (\ref{inehu}). According to
(\ref{hltgv}), the right-hand side of (\ref{grofu}) and
(\ref{cuphu}) intersect $T^\bbC v$ by the set
\begin{eqnarray*}
T\exp(P_X)v=\exp(Q_X)T\exp(C_T)v;
\end{eqnarray*}
clearly, it is dense in (\ref{cuphu}). Since $Q_X$ is compact, the
set (\ref{cuphu}) is closed. The compactness of $Q_X$, the above
equality, and Lemma~\ref{cloct} imply that (\ref{addfa}) is
closed; hence, it is the same as the right-hand side of
(\ref{grofu}).

Each of the sets (\ref{grofu})--(\ref{addfa}) includes $Y$. Thus,
it remains to prove that (\ref{addfa}) is polynomially convex. If
$x\in\frt^\bbR\setminus P_X$, then there exists $s\in\bbR^n$ such
that
\begin{eqnarray}\label{sepa}
\sup\{\scal{y}{s}:\,y\in P_X\}<\scal{x}{s}.
\end{eqnarray}
Since $Q_X$ is compact, the linear functional on $\frt$ in the
right-hand side of (\ref{sepa}) must be nonnegative on $C_T$.
According to (\ref{defct}), we may assume that
$s\in\clos\bbR^n_+$. It follows that (\ref{sepa}) holds in a
neighborhood of $s$ in $\clos\bbR^n_+$. Thus, $s$ can be assumed
rational (hence, integer) with strictly positive entries. Then,
$p(z)=z_1^{s_1}\dots z_n^{s_n}$ is a holomorphic polynomial such
that $|p\,|$ separates $\exp(x)v$ and $T\clos(\exp(P_X)v)$.
Therefore,
\begin{eqnarray*}
\wh Y\cap T^\bbC v=T\exp(P_X)v.
\end{eqnarray*}
For any $\iota \in I_T$, the projection $z\to \iota z$ commutes
with $T$. This makes it possible to apply the above arguments to
the vector $\iota v$, the set $\iota X$, and to the restriction of
$T$ to $\iota \bbC^n$. Consequently,
\begin{eqnarray}\label{iotor}
\wh{\iota Y}\cap T^\bbC \iota v=T\exp(P_{\iota X})\iota v = \iota
T\exp(P_{X})v
\end{eqnarray}
(clearly, $\iota\exp(P_{X}) v=\exp(P_{\iota X})\iota v$). By
(\ref{defit}), $\iota Y\subseteq\wh Y$, hence, $\wh{\iota
Y}\subseteq\wh Y$; on the other hand, $\iota \wh
Y\subseteq\wh{\iota Y}$ since $p\circ\iota$ is a polynomial on
$\bbC^n$ for any polynomial $p\,$ on $\iota\bbC^n$. Thus, $\iota
\wh Y=\wh{\iota Y}=\wh Y\cap\iota\bbC^n$. Together with
(\ref{iotor}), this implies the polynomial convexity of
(\ref{addfa}).
\end{proof}
If $T=\bbT^n$, then Proposition~\ref{hulfi} follows from the
well-known characterization of polynomially convex Reinhardt
domains.
\begin{corollary}\label{noclo}
For any $v\in\csn$, the orbit $T^\bbC v$ is closed in $\bbC^n$ if
and only if $\,Tv$ is polynomially convex, and this is equivalent
to $C_T=0$. Then, $\wh{Y}=T\exp(Q_X)v$ for all $Y,X$ as above.
\end{corollary}
\begin{proof}
The orbit  $T^\bbC v$ is closed if and only if the convex hull of
$\spec(v)=\spec(\bbC^n)$ contains $0$ in its relative interior
(see, for example, \cite[Proposition~6.15]{VP}).  Since
$T\subseteq\GL(n,\bbC)$, the set $\spec(\bbC^n)$ is generating in
$\frt^*$. Hence, $T^\bbC v$ is closed if and only if $C_T=0$; by
(\ref{hultv}),  this is equivalent to $\wh{Tv}=Tv$. Then,
$\wh{Y}=T\exp(Q_X)v$ by (\ref{grofu}) and (\ref{speceq}).
\end{proof}
There is a version of the first assertion for an arbitrary compact
linear group $G$: a $G^\bbC$-orbit is closed if and only if it
contains a polynomially convex $G$-orbit (\cite[Theorem~1 and
Theorem~5]{GL}). For a torus $T$, all $T$-orbits in $T^\bbC v$ are
simultaneously polynomially convex or non-convex, but this is not
true if $G$ is not abelian.

\section{Finite extensions of $T$ that keep a $T^\bbC$-orbit}
In this section, we consider the case where the set $X$ defined in
the previous section is an orbit of a finite group $F$ which
normalizes $T$ and keeps the $T^\bbC$-orbit. We assume that
$T\subseteq G$, $T$ is a torus, $G$ is a subgroup of $\GL(V)$, $F$
is a finite subgroup of $G$, and
\begin{eqnarray}
&G=FT=TF,\quad F\cong G/T,\label{injec}\\
&Gv\subseteq T^\bbC v,\label{incor}\\
&v\in\csn\subset\bbC^n=V.\label{vgene}
\end{eqnarray}
By (\ref{injec}), $T$ is normal in $G$. Clearly, (\ref{incor}) is
equivalent to $Fv\subseteq T^\bbC v$ and to the connectedness of
$G^\bbC$. Here is an illustrating example.
\begin{example}\label{autd2}
Let $G=\Aut_0(\bbD^2)$ be the group of linear automorphisms of the
bidisc $\bbD^2\subset\bbC^2$. Clearly,  $G=FT$, where $F=S_2$ is
generated by the transposition $\tau$ of the coordinates,
$T=\bbT^2$, $T^\bbC=(\bbC^*)^2$, and $T^\bbC v=(\bbC^*)^2$ for any
$v$ that lies outside the coordinate lines. Thus, (\ref{incor})
holds for all $v\in(\bbC^*)^2$ (however, (\ref{incor}) fails for
any $v\neq0$ in $\bbC^2\setminus (\bbC^*)^2$). The hull $\wh{Gv}$
can be distinguished by the inequalities
\begin{eqnarray}
&\max\{|z_1|,|z_2|\}\leq \max\{|v_1|,|v_2|\},\label{ined2}\\
&|z_1z_2|\leq|v_1v_2|.\label{ined2a}
\end{eqnarray}
Clearly, (\ref{ined2}) and (\ref{ined2a}) define a polynomially
convex set. Let $z_1,z_2>0$ (a generic $T$-orbit evidently
contains such a point $z$). Then, $z$ and $\tau z$ can be joined
by an analytic strip with the boundary in $Tz\cup T\tau z$:
\begin{eqnarray*}
\la_z(s)=(z_1^{1-s}z_2^{s},z_1^{s}z_2^{1-s}), \quad s\in S.
\end{eqnarray*}
Set $q=\ln\frac{z_1}{z_2}$ and let $z_1>z_2$. Then, the strip can
be written in the form
\begin{eqnarray*}
\la_z(s)=(e^{-s}z_1,e^{s}z_2),\quad 0\leq\re s\leq q.
\end{eqnarray*}
It is periodic with the period $2\pi i$ and defines a
$\tau$-invariant annulus in $\wh{Gv}$ with $\tau$-fixed points
$(\sqrt{z_1z_2},\sqrt{z_1z_2})$ and
$(-\sqrt{z_1z_2},-\sqrt{z_1z_2})$. As $z_2\to0$, the annulus tends
to a couple of discs: $(e^{-s}z_1,0)$ and $(0,e^{-s}z_1)$, where
$\re s>0$, $0\leq\im s\leq2\pi$ (the circle $\re s=\frac{q}{2}$,
$0\leq\im s\leq 2\pi$ collapses to zero). Let $z\in
\wh{Gv}\cap\bbR^2$. Then $\wh{Gv}$ contains a bidisc $\bbD^2z$. It
intersects $\bbR^2$ by a rectangle, which is symmetric with
respect to the coordinate axes. If $z$ lies on an axis, then the
rectangle degenerates into a segment. Let $v_1> v_2>0$. The union
of these rectangles with vertices in the set $Q$ of real points of
the annulus, which joins $v$ and $\tau v$, is a curvilinear
octagon. It degenerates into a pair of segments if $v_2=0$ and
into a square if $v_1=v_2$  (see \cite[Fig. 2]{KZ} for the
3-dimensional case). In the logarithmic coordinates in the first
quadrant, $Q$ is a segment. Also, note that all nontrivial
$T^\bbC$-orbits are not closed.\qed
\end{example}
In \cite{Bj}, Bj\"ork found a typical situation where analytic
annuli appear in the maximal ideal space $\cM_A$ of a commutative
Banach algebra $A$ which admits a nontrivial action of $\bbT$ by
automorphisms: this happens if $T$-invariant functions on $\cM_A$
do not separate distinct $\bbT$-orbits. In \cite{GL}, it was noted
that analytic strips and/or annuli appear in $\wh{Gv}$ if the
stable subgroup of $v$ in $G^\bbC$ does not coincide with the
complexification of the stable subgroup of $v$ in $G$.
\begin{proposition}\label{autdn}
The hulls $\wh{Gv}$ for orbits of $G=\Aut_0(\bbD^n)=S_n\bbT^n$ are
distinguished by inequalities {\rm(\ref{mukin})}, where $\mu_k$
are defined by {\rm(\ref{mukde})}.
\end{proposition}
\begin{proof}
The approximation by decreasing sequences of hulls makes it
possible to reduce the proposition to the case of a generic $v$ in
{\rm(\ref{mukin})}. Then, applying to $v=(v_1,\dots,v_n)$ a
suitable transformation in $\bbT^n$, we may assume that
\begin{eqnarray}\label{gewey}
v_1> v_2> \dots> v_n> 0.
\end{eqnarray}
Moreover, we may use Proposition~\ref{hulfi} with $X=S_n v$,
$C_T=-\clos\bbR^n_+$ (we keep the notation of
Proposition~\ref{hulfi}). Since $X$, $Q_X$, $C_T$, $P_X$, and
$\mu_k$ are $S_n$-invariant, $S_n$ is transitive on $X$, by
(\ref{mukde}), (\ref{mukin}), and (\ref{inehu}), it is sufficient
to prove that the vectors $\xi_k=\sum_{r=1}^k\eps_r$,
$k=1,\dots,n$, correspond to the faces of $P_X$ that meet at $v$,
are orthogonal to them, and look outside of $P_X$.

Set $\eta_1=\eps_2-\eps_1,
\dots,\eta_{n-1}=\eps_{n}-\eps_{n-1},\eta_n=-\eps_n$. Then
$\{-\eta_k\}_{k=1}^n$ is a base in $\bbR_n$, which is dual to the
base  $\{\xi_k\}_{k=1}^n$. We claim that the cone of the
polyhedron $P_X$ at the vertex $v$ is generated by
$\{\eta_k\}_{k=1}^n$. This implies the assertion above (note that
both cones are simplicial). If $\tau\in S_n$ is a transposition
$(k,j)$, then $v-\tau v=(v_k-v_j)(\eps_k-\eps_j)$. If $\si,\ka\in
S_n$ then $v-\si\ka v=(v-\ka v)+(\ka v-\si\ka v)$. Furthermore,
$S_n$ is generated by transpositions $(k,k+1)$, and
$v_k-v_{k+1}>0$ by (\ref{gewey}), where $k=1,\dots,n-1$.
Therefore, vectors $\eta_1,\dots,\eta_{n-1}$ generate the cone of
$Q_X$ at $v$. Since $-\eps_k=\sum_{r=0}^{k-1}\eta_{n-r}$ and $C_T$
is generated by $-\eps_k$, $k=1,\dots,n$, this proves the
proposition.
\end{proof}

Property (\ref{incor}) implies $\spann_\bbC Gv=\spann_\bbC T^\bbC
v$. Hence, we may assume that (\ref{sispe}) is valid. Then, a
generic $\xi\in\frt$ has a simple spectrum. Any $f\in F$ permutes
eigenvalues and eigenspaces. Thus, assuming (\ref{sispe}) and
identifying $V$ with $C(K)$ in accordance with (\ref{vcnck}), we
get that each element of $F$ is a composition of a permutation of
$K$ and a multiplication by a function on $K$. Further,
(\ref{vgene}) implies that the stable subgroup of $v$ in $T^\bbC$
is trivial. Hence,
\begin{eqnarray*}
T^\bbR v\cong T^\bbC/T\cong \frt^\bbR,
\end{eqnarray*}
where the identification of $T^\bbR v$ and $\frt^\bbR$ is realized
by $\xi\to\exp(\xi)v$, $\xi\in\frt^\bbR$.
\begin{lemma}\label{fixtc}
Let $G\subset\GL(V)$, a subgroup $F\subseteq G$, and $v\in V$
satisfy (\ref{injec})--(\ref{vgene}). Then $T^\bbC v$ contains a
$G$-invariant $T$-orbit. Moreover, there exists a mapping $f\to
t_f$, $F\to T$, such that $\wt F=\{t_f f:\,f\in F\}$ is a subgroup
of $G$ which has a fixed point in $T^\bbC$ and satisfies
(\ref{injec})--(\ref{vgene}).
\end{lemma}
\begin{proof}
The group $F$ naturally acts on $T^\bbR v\cong
T^\bbC/T\cong\frt^\bbR$. Any $g\in F$ is a composition of $\si\in
S_K$ and a multiplication by a function in $C(K)$. Since $\frt$
acts on $C(K)$ by multiplication on linear functions and $\si$
induces a linear transformation in $\frt^\bbC$, the induced action
of $F$ on $\frt^\bbR$ is affine. Since $F$ is finite, it has a
fixed point in $\frt^\bbR$. Hence, $T^\bbC v$ contains a
$G$-invariant $T$-orbit. Let us fix a point $u$ in it and define
$t_f$ by $t_ffu=u$; the choice is unique due to (\ref{vgene}).
Taken together with (\ref{injec}), this implies that $\wt F$ is a
group, which obviously satisfies the lemma.
\end{proof}
According to Lemma~\ref{fixtc}, we may assume without loss of
generality that
\begin{eqnarray}\label{fixev}
fv=v\quad\mbox{for all}\quad f\in F.
\end{eqnarray}
In the following example we give a construction (associated with a
given finite group $F$) for orbits with property (\ref{incor}).

\begin{example}\label{maico}
Let $\frt$ be a real linear space, $\frt^*$ be the dual space to
$\frt$, $L$ be a lattice in $\frt$, and $L^*\subset\frt^*$ be the
dual lattice to $L$. Set
\begin{eqnarray*}
\la_x(y)=y(x),\quad\mbox{where}\ x\in\frt,\ y\in\frt^*.
\end{eqnarray*}
Let $K$ be a finite subset of $L^*$ that generates $L^*$ as a
subgroup of the vector group $\frt^*$. Then
\begin{eqnarray}
&\frt^*=\spann_\bbR K,\label{spank}\\
&L=\{x\in\frt:\,\la_x(K)\subset\bbZ\}.\label{defz}
\end{eqnarray}
Further, let $F$ be a finite subgroup of $\GL(\frt)$ which keeps
$K$. Set $V=C(K)$. The mapping
\begin{eqnarray}\label{defla}
\la:\,x\to i\la_x\big|_K
\end{eqnarray}
is an embedding $\frt\to V$, which has a natural extension to
$\frt^\bbC$. Set
\begin{eqnarray}\label{defex}
\exp(x)=e^{2\pi i\la_x}.
\end{eqnarray}
Clearly, $L=\ker\exp$. Hence, $\exp$ defines an embedding of
$T=\frt/L$ and $T^\bbC$ into the group $\csn$:
\begin{eqnarray*}
T^\bbC=\exp(\frt^\bbC)\subseteq \csn.
\end{eqnarray*}
The group $T^\bbC$ acts on $C(K)$ by multiplication. The inclusion
$v\in C(K)^F$ is the same as (\ref{fixev}); it implies
(\ref{incor}). Furthermore, if $v\in \csn$, then
\begin{eqnarray}\label{spatc}
\spann_\bbC Tv=V.
\end{eqnarray}
Indeed, the space $\spann_\bbC T$ is a subalgebra of $C(K)$, which
separates points of the finite set $K$.  Hence, it coincides with
$C(K)$. \qed
\end{example}
\begin{theorem}\label{isase}
Let a group $G\subset\GL(V)$, a finite subgroup $F\subseteq G$, a
torus $T$, and a vector $v\in V$ satisfy
(\ref{injec})--(\ref{vgene}), (\ref{fixev}), and (\ref{spatc}).
Then $V,G,F,T,v$ can be realized as in Example~\ref{maico}, where
\begin{eqnarray}\label{vgenf}
v\in \csn\cap C(K)^F.
\end{eqnarray}
Conversely, if $V,G,F,T,v$  are as in Example~\ref{maico} and $v$
satisfies (\ref{vgenf}), then (\ref{injec})--(\ref{vgene}),
(\ref{fixev}), and (\ref{spatc}) are true.
\end{theorem}
\begin{proof}
The group $F$ acts in $\frt$ and $\frt^*$ by the adjoint action.
Let $K\subset\frt^*$ be the collection of all weights for the
representation of $T$ in $V$; clearly, $K$ is $F$-invariant. It
follows from (\ref{spatc}) and (\ref{incor}) that the weights are
multiplicity free. This defines an equivariant linear isomorphism
between $V$ and $C(K)$, where the group $T$ acts by
multiplication. Thus, $\la$ and $\exp$ are well defined by
(\ref{defla}) and (\ref{defex}). According to (\ref{fixev}) and
(\ref{vgene}), (\ref{vgenf}) is true; (\ref{spank}) holds since
$T\subset\GL(V)$ is compact and acts effectively on $V$ (note that
the annihilator of $\spann_\bbR K$ in $\frt$ acts trivially due to
(\ref{defex}) and (\ref{vgenf})). Let us define $L$ by
(\ref{defz}). Then $L=\ker\exp$ by (\ref{defex}). Hence, $L$ is a
lattice in $\frt$ and the group $L^*$ generated by $K$ is the dual
lattice in $\frt^*$.

The converse was proved in Example~\ref{maico}.
\end{proof}

\section{Finite extensions of $T$ which keep generic $T^\bbC$-orbits}

In what follows, we use the setting of Example~\ref{maico}. Let
$Z$ denote the centralizer of $G$ in $\GL(V)$. We assume that
$\csn$ acts in $V=C(K)$ by multiplication.
\begin{lemma}\label{centr}
$Z=C(K)^F\cap \csn$.
\end{lemma}
\begin{proof}
Since $\la(\frt)$ separates points of $K$, $Z\subseteq\csn$. 
The multiplication by $u\in C(K)$ commutes with $F$ if and only if
$u$ is $F$-invariant.
\end{proof}

In general, condition (\ref{incor}) does not hold for a generic
vector $v$. Hence, there is a natural problem: {\sl describe $V$
and $G$ such that generic orbits satisfy (\ref{incor})}. The
following proposition contains a simple criterion.
\begin{proposition}\label{sicri}
Let $V,G$ be as in Example~\ref{maico}. Then $G$ satisfies
(\ref{incor}) for a generic $v\in V$ if and only if
\begin{eqnarray}\label{opfix}
C(K)=\la(\frt^\bbC)+C(K)^F.
\end{eqnarray}
In this case, each $T^\bbC$-orbit in $\csn$ intersects $C(K)^F$.
\end{proposition}
\begin{proof}
It follows from Lemma~\ref{centr} that the right-hand side of
(\ref{opfix}) is the tangent space at $\one$ to the set
$T^\bbC\,Z$. Clearly, $\wt G^\bbC=ZG^\bbC$ is a group, $T^\bbC Z$
is the identity component of $\wt G^\bbC$, and the right-hand side
of (\ref{opfix}) is the tangent space to $\wt G^\bbC\one$. Hence,
(\ref{opfix}) holds if and only if $\wt G^\bbC\one$ is open.
Moreover, this is equivalent to the equality
$T^\bbC Z=\exp(\la(\frt^\bbC)+C(K)^F)=\csn$. 
Therefore, each $T^\bbC$-orbit in $\csn$ intersects $C(K)^F$,
i.e., contains an $F$-fixed point. Thus, (\ref{opfix}) implies
(\ref{incor}) for $v\in \csn$.

Let (\ref{incor}) hold and let $W$ be an $F$-invariant
neighborhood of $\one$. If $W$ is sufficiently small, then the
condition $\log\one=0$ defines a branch of $\log$ in $W$. We may
assume that $\log W$ is convex and symmetric.
This makes it possible to define roots in $W$:
$w^{\frac1r}=\exp\left(\frac1r\log w\right)$.
For $v\in W^{\frac12}$ and $f\in F$, set $g_f=\left(\frac
{fv}{v}\right)^{\frac1r}$, where $r=\card F$, and
$g=\prod\nolimits_{f\in F}g_f$.
Then $g v$ is $F$-fixed. If (\ref{incor}) holds for $v$, then
$g_f\in T^\bbC$ for all $f\in F$; hence, $g v\in T^\bbC v$.
Consequently, for all $v\in W$,  $T^\bbC v$ intersects $C(K)^F$.
Since $Z$ keeps this property of orbits, it follows that $T^\bbC
Z$ has a nonempty interior. This implies (\ref{opfix}).
\end{proof}
\begin{theorem}\label{charf}
Let $G\subset\GL(V)$ be a semidirect product of a torus $T$ and a
finite subgroup $F$, and let $Z$ be the centralizer of $G$ in
$\GL(V)$. Suppose that $\spann_\bbC Tv=V$ for some $v\in V$. Then
the following conditions are equivalent:
\begin{itemize}
\item[(\romannumeral1)] $Gv\subset T^\bbC v$ for a generic $v\in
V$; \item[(\romannumeral2)] $G^\bbC Zv$ is open in $V$  for a
generic $v\in V$.
\end{itemize}
\end{theorem}
\begin{proof}
By Theorem~\ref{isase}, we may use the construction of
Example~\ref{maico}. According to Lemma~\ref{centr},
(\romannumeral2) is equivalent to
(\ref{opfix}), and the assertion follows from
Proposition~\ref{sicri}.
\end{proof}
We shall give a constructive description of these spaces and
groups. Set
\begin{eqnarray*}
C_0(K)=\left\{u\in C(K):\,\ \sum\nolimits_{q\in K}u(q)=0\right\}.
\end{eqnarray*}
Sometimes, we identify points in $K$ with their characteristic
functions.
\begin{example}\label{maict}
Let $V=\bbC^n=C(K)$, where $K=\{1,\dots,n\}$, let $F$ be a
subgroup of $S_n$, and
\begin{eqnarray}\label{deck}
K=K_1\cup\dots\cup K_p
\end{eqnarray}
be the partition of $K$ into $F$-orbits. For $k\in\{1,\dots,p\}$,
set $V_k=C(K_k)$. Then $V=V_1\oplus\dots\oplus V_p$.  Set
\begin{eqnarray*}
&\frt_k^{0}=C_0(K_k)\cap i\bbR^n,\\
&T_k^0=\exp(\frt_k^0)\subset C(K_k),
\end{eqnarray*}
where $\exp$ is defined by (\ref{defex}). Set
$\frt^0=\frt_1^0\oplus\dots\oplus \frt_p^0,$
\begin{eqnarray*}
T^0=\exp(\frt^0)=T_1^0\times\dots\times T_p^0.
\end{eqnarray*}
Let $T$ be an $F$-invariant torus such that
\begin{eqnarray}\label{inct}
T^0\subseteq T\subseteq\bbT^n
\end{eqnarray}
and set $G=FT$. Then generic $G^\bbC$-orbits satisfy
(\ref{incor}). The group $G$ is irreducible if and only if $F$ is
transitive on $K$; in general, $F$-orbits in $K$ define
$G$-irreducible components of $V$.  There are two extreme cases in
(\ref{inct}).
\begin{itemize}
\item[(A)]\label{exama} If $T=\bbT^n$, then there is one open
orbit $\csn$ of the group $G^\bbC=FT^\bbC$, which evidently
satisfies (\ref{incor}). If $F$ is nontrivial, then there exist
degenerate orbits that do not satisfy (\ref{incor}); moreover, if
$F$ is transitive on $K$, then all non-open $G^\bbC$-orbits,
except for zero, are nontrivial finite unions of $T^\bbC$-orbits.
\item[(B)]\label{examb} If $T=T^0$, then generic orbits are
closed. They have codimension $p$ and are distinguished by
equations
\begin{eqnarray*}
\prod\nolimits_{r\in K_k}z_r=c_k,
\end{eqnarray*}
where $c_k\in\bbC^*$, $k=1,\dots,p$.
\end{itemize}
\end{example}
Note that (A) and (B) are invariant under the Cartesian product
(the group $F$ need not be the product of groups $F_k$ of
irreducible components but must have the same orbits in $K$ as
$F_1\times\dots\times F_p$). In terms of Example~\ref{maico}: in
(A), $\frt=\bbR^n$, the mapping $\la:\,\frt^\bbC\to C(K)$ is
surjective, $K=\{\eps_1,\dots,\eps_n\}$; in (B), $\frt=i\bbR^n\cap
C_0(K)$, $\la(\frt^\bbC)=C_0(K)$,  and the set $K$ is the
projection of $\{\eps_1,\dots,\eps_n\}$ into $\frt^*=\frt$. In
both cases, $K$ is the set of all vertices of a regular simplex.
\qed
\begin{theorem}\label{asexa}
Let $V,G$ be as in Theorem~\ref{charf} and let
{\rm{(\romannumeral1)}} hold. Then $V,G$ can be realized as in
Example~\ref{maict}. Furthermore,
\begin{itemize}
\item[\rm(1)] $V,G$ are of type {\rm(A)} if and only if $G^\bbC$
has an open orbit, \item[\rm(2)] {\rm(B)} is equivalent to the
assumption that the center of $G$ is finite, \item[\rm(3)] if $G$
is irreducible, then either {\rm(A)} or {\rm(B)} holds.
\end{itemize}
\end{theorem}
Let $C(K)^F_+$ be the cone of all nonnegative functions in
$C(K)^F$.
\begin{lemma}\label{crifi}
Let $G$ and $V$ be as in Example~\ref{maico}. Then, the orbit
$G^\bbC v$ is closed for a generic $v\in V$ if and only if
\begin{eqnarray}\label{close}
\frt^\bbR\cap C(K)^F_+=0.
\end{eqnarray}
\end{lemma}
\begin{proof}
Clearly, $G^\bbC v$ is closed if and only if $T^\bbC v$ is closed.
Let $v\in\csn$.  By Proposition~\ref{hulfi} and
Corollary~\ref{noclo}, $T^\bbC v$ is not closed if and only if
$C_T\neq0$. Since $C_T$ is $F$-invariant by (\ref{defct}), it
contains $\sum_{f\in F} fu$ for each $u\in C_T$. Thus, $C_T=0$ is
equivalent to (\ref{close}).
\end{proof}
\begin{proof}[Proof of Theorem~\ref{asexa}]
Suppose that $G$ is irreducible or, equivalently, $F$ is
transitive. Then $Z=\bbC^*\,\one$ according to Lemma~\ref{centr}.
If $\one\in\la(\frt^\bbC)$, then $T^\bbC\supseteq Z$ and $T^\bbC
v$ is open for a generic $v\in V$ by Theorem~\ref{charf}. If
$\one\notin\la(\frt^\bbC)$, then (\ref{close}) is true; by
Lemma~\ref{crifi}, $T^\bbC v$ is closed for a generic $v\in V$. By
Proposition~\ref{sicri}, a generic $T^\bbC$-orbit intersects
$C^*\one$. Consequently, we have
\begin{eqnarray}\label{codim}
\codim G^\bbC v=1.
\end{eqnarray}
Let $\one\in T^\bbC\cap Z$. The orthogonal projection of $\one$
into the tangent space ${\mathrm T}_{\one}T^\bbC\one$ is
$F$-fixed. Hence, it is proportional to $\one$; since
$\one\notin\la(\frt^\bbC)$, this implies
$\one\perp\la(\frt^\bbC)\one$. Therefore, ${\mathrm T}_\one
T^\bbC\one$ coincides with the tangent space to the hypersurface
$z_1\dots z_n=1$ at $\one$;  since the monomial on the left is an
eigenfunction of $T^\bbC$, this group keeps it. Due to
(\ref{codim}), $T^\bbC\one$ coincides with this hypersurface.
Then, $T=\bbT^n\cap\SU(n)$, and any $T^\bbC$-orbit that intersects
$Z$ is a hypersurface $z_1\dots z_n=c$, for some $c\in\bbC^*$.
This implies $\frt^\bbC=C_0(K)$ and $T=T^0$.

Thus, the theorem is proved for all irreducible $G$. The
projection onto each irreducible component keeps the property
(\ref{incor}) for generic orbits since it commutes with $G$.
Hence, (\romannumeral1) holds for all irreducible components. They
correspond to $F$-orbits $K_k$ in the partition (\ref{deck}).
Let $\frt_k^0$, $k=1,\dots,p\,$, be defined as in
Example~\ref{maict}. According to the arguments above,
$\la\left(\frt|_{K_k}\right)\supseteq\la(\frt_k^0)$ for all $k$.
If $x\in\frt$, then the averaging
\begin{eqnarray*}
&Ax=\frac1r\sum\nolimits_{f\in F}fx,\quad r=\card F,
\end{eqnarray*}
distinguishes the $F$-fixed component of $x$ (i.e., $Ax\in
C(K)^F\cap\frt$ and $x-Ax\in\frt^0$); since $\frt$ is
$F$-invariant, it contains both components. By Lemma~\ref{centr},
if $G$ has a finite center, then
$\la\left(\frt|_{K_k}\right)=\la(\frt_k^0)$ for all $k$. It
follows that
$$\frt\subseteq\frt^0=\frt_1^0\oplus\dots\oplus \frt_p^0.$$
On the other hand, (\romannumeral2) and Lemma~\ref{centr} imply
$\codim \frt\leq\dim C(K)^F=p$. Hence, the inclusion above is in
fact the equality. Thus, we get (B) assuming that $G$ has a finite
center. The converse is true since $\frt^0$ does not contain a
nontrivial $F$-fixed element. The same arguments show that any
$F$-invariant torus $T$ includes $T^0$ if (\romannumeral1) is
true. This proves that $V,G$ admit the realization of
Example~\ref{maict}; (1) and (2) are clear.
\end{proof}
\begin{corollary}
Let $G$ be as in Theorems~4.5 and 4.3. Then $G$ contains a closed
subgroup $G^0$ such that
\begin{itemize}
\item[(1)] each connected component of $G$ contains a connected
component of $G^0$, \item[(2)] $G^0$ has a finite center,
\item[(3)] generic orbits of $(G^0)^\bbC$ are closed, \item[(4)]
$Gv\cap T^\bbR v=G^0v\cap(T^0)^\bbR v$ for a generic $v\in V$.
\end{itemize}
\end{corollary}
\begin{proof}
By Theorem~\ref{asexa} and (\ref{inct}), $G\supseteq T^0$, where
$T^0$ is as in (B). Clearly, $F$ normalizes $T^0$. Hence,
$G^0=FT^0$ is a group, which satisfies the corollary.
\end{proof}
Proposition~\ref{hulfi} makes it possible to find $\wh{Gv}$ for
$G$ as above.  If $T=\bbT^n$, then $T\supset\bbZ_2^n$ and generic
$T$-orbits intersect $\bbR^n_+$; hence, we may assume
$v\in\bbR^n_+$. Then $\wh{Tv}\cap\bbR^n$ is a parallelepiped
$\Pi_v=\conv\{(\pm v_1,\dots,\pm v_n)\}$. Clearly,
$\Pi_v=\bbZ^n_2\Pi_v^+$, where $\Pi_v^+=\Pi_v\cap\clos\bbR^n_+$.
Since $\bbR^n_+=T^\bbR v$, we may use Proposition~\ref{hulfi} with
$X=Fv$, $C_T=-\clos\bbR^n_+$, $P_X=\conv(Fv)-\bbR^n_+$:
\begin{eqnarray*}
\wh{Gv}=\cup_{u\in\exp(Q_v)}\bbD^n u=T\cup_{u\in\exp(Q_v)}\Pi_{u}=
T\cup_{u\in\exp(Q_v)}\Pi_{u}^+,
\end{eqnarray*}
where $Q_v=\conv{Fv}$. For the description in the form
(\ref{inehu}), one has to know normal vectors to faces of $\conv
Fv$. Since $F$ may be an arbitrary subgroup of $S_n$, they need
not be proportional to rational vectors (for example, this is true
for the cyclic subgroup of order 3 in $S_3$). We shall describe
the situation where they are locally independent of $v$; since
they depend on $v$ continuously, this is equivalent to the
condition that they are rational. Note that the vector which joins
two points in $\frt^\bbR$ as in Lemma~\ref{strih} is rational if
and only if the strip reduces to an annulus. In
Example~\ref{maict}, $F$ need not be the product of groups
corresponding to the irreducible components; we shall see that $F$
possesses this property in the case under consideration.

Let $U$ be a real vector space and $F\subset\GL(U)$ be a finite
group. Set
\begin{eqnarray*}
C_u=\cone(u-Fu);
\end{eqnarray*}
this is the cone at the vertex $u$ of the polytope $\conv(Fu)$
(which may be degenerate). We say that {\it $C_u$ is locally
independent of $u$} if, for a generic $u\in U$, $C_u=C_w$ for all
$w$ that are sufficiently close to $u$.
\begin{lemma}\label{coxno}
Let $U$ be a real vector space and $F$ be a finite subgroup of
$\GL(U)$. Suppose that $C_u$ is locally independent of $u$. Then
$F$ is generated by reflections in hyperplanes in $U$.
\end{lemma}
\begin{proof}
We may assume without loss of generality that $U$ is equipped with
an inner product and that $F\subseteq\OO(U)$. Let $\bbR_+(u-fu)$,
$f\in F$, be an extreme ray of $C_u$. The equality $C_u=C_w$ for
$w$ in a neighborhood of $u$ implies that this ray does not change
near $u$. Hence, $\dim(\one-f)U=1$. Since $f$ is orthogonal and
nontrivial, it is a refection in a hyperplane. The stable subgroup
of a generic $u\in U$ is trivial (hence, $F$ acts freely on a
generic orbit) and each vertex of $\conv(Fu)$ can be joined with
$u$ by a chain of edges. Applying the above arguments repeatedly
to $u,fu$, etc., we get that $F$ is generated by reflections in
hyperplanes.
\end{proof}
For any $g\in\bbZ_2^nS_n$ and $k=1,\dots,n$,
$g\eps_k=\pm\,\eps_{\si(k)}$ for some $\si\in S_n$. The mapping
$f\to\si$ is a natural homomorphism $\bbZ_2^nS_n\to S_n$, which we
denote by $\phi$.
\begin{lemma}\label{perre}
Let $F$ be a transitive subgroup of $S_n$ acting in $\bbR^n$ by
permutations of coordinates and let a group $H\subseteq
\bbZ_2^nS_n$ be generated by reflections in hyperplanes in
$\bbR^n$. If $\phi(H)=F$, then $F=S_n$.
\end{lemma}
\begin{proof}
Let $\rho$ be a reflection in a hyperplane in $\bbR^n$. If
$\rho\in\bbZ_2^nS_n=BC_n$, then it is conjugate to a reflection in
a wall of the Weyl chamber that is distinguished by the
inequalities $x_1>\dots>x_n>0$. Hence, $\phi(\rho)$ is a
transposition if it is nontrivial. Since $F=\phi(H)$, $F$ is
generated by transpositions. It remains to note that any subgroup
of $S_n$, which is generated by transpositions, coincides with
$S_n$ if it is transitive on $\{1,\dots,n\}$ (consider the graph
with the vertices $\{1,\dots,n\}$ and edges corresponding to
transpositions and note that inclusions $(k,l)\in F$, $(l,m)\in F$
imply $(k,m)\in F$; this makes it possible to use the induction).
\end{proof}
We say that a pair $(V,G)$ is {\it standard} if it is isomorphic
to (A) or (B) in Example~\ref{maict} with $F=S_K$. The {\it
product} of pairs $(V_k,G_k)$, $k=1,\dots,m$, is the pair
$(\sum_{k=1}^m V_k, \prod_{k=1}^m G_k)$.
\begin{theorem}\label{huini}
Let $G=FT$ be a compact subgroup of $\GL(n,\bbC)$, where
$T\subseteq\bbT^n$ is a torus and $F$ is a subgroup of $S_n$.
Suppose that $Gv\subset T^\bbC v$ for a generic $v\in V$ and
\begin{itemize}
\item[\rm(1)] either $T=\bbT^n$ or the center of $G$ is finite,
\item[\rm(2)] for a generic $v\in\bbC^n$,  $\wh{Gv}$ can be
distinguished in $\clos T^\bbC v$ by a family of inequalities
\begin{eqnarray*}
|z_1|^{s_1}\dots |z_n|^{s_n}\leq \rho_s(v),
\end{eqnarray*}
where $\rho_s(v)\geq0$ and vector $s=(s_1,\dots,s_n)$ runs over a
certain finite subset of $\bbR^n$ which is independent of $v$.
\end{itemize}
Then $(V,G)$ is isomorphic to the product of standard pairs.
Moreover, if $G$ is irreducible, then $(V,G)$ is standard.
\end{theorem}
\begin{proof}
Let $G$ be irreducible. Then $F$ is transitive and  $(V,G)$ are as
in (A) or as in (B) by Theorem~\ref{asexa}.  Suppose that (B) is
the case. It follows from (2) and Proposition~\ref{hulfi} that the
polytope $Q_X\subset\frt^\bbR$, where $X=Gv\cap T^\bbR v$, for a
generic $v$, satisfies the assumption of Lemma~\ref{coxno}.
Therefore, $F$ is generated by reflections (we may assume that
$F\subset\OO(\frt^\bbR)$). They extend to reflections in
hyperplanes in $\frt^\bbR+\bbR\,\one=\bbR^n$ if we assume that
they fix $\one$. Then, Lemma~\ref{perre} implies $F=S_n$. The case
(A) can be reduced to (B): it is sufficient to replace $\bbT^n$
with $T=\SU(n)\cap\bbT^n$ since $F$ evidently keeps $T$ and to
note that (2) remains true due to Proposition~\ref{hulfi}. Thus,
$(V,G)$ is standard.

Let the center of $G$ be finite. According to Theorem~\ref{asexa},
$T$ may be identified with the group $T^0$ in Example~\ref{maict}.
In particular, $G^\bbC v$ is closed for a generic $v$ and $C_T=0$
due to Proposition~\ref{hulfi}. By Proposition~\ref{sicri},
generic orbits contain $F$-fixed points. Applying the arguments
above (which did not use the assumption that $G$ is irreducible),
we get that the cones at the vertices of the convex polytope
$Q_X$, $X=Gv\cap T^\bbR\subset\frt^\bbR$, are locally independent
of $v$. Clearly, the same is true for its projection into each
space $\frt^0_k$ corresponding to an irreducible component $V_k$
of $V=\bbC^n$. This implies that all irreducible components are
standard. Thus, $F_k=S(K_k)$, where $k=1,\dots,p$ and
$K=K_1\cup\dots\cup K_p$ is the partition of $K$ into $F$-orbits.
Due to Theorem~\ref{asexa}, it is sufficient to prove that
\begin{eqnarray}\label{produ}
F=F_1\times\dots\times F_p.
\end{eqnarray}
By Lemma~\ref{coxno}, $F|_{\frt^0}$ is generated by reflections in
hyperplanes in $\frt^0$; the condition that they keep real
$F$-invariant functions on $K$ uniquely defines their extension to
$\bbR^n$. Hence, $F$ is generated by reflections in $\bbR^n$. A
permutation which induces a reflection in a hyperplane in $\bbR^n$
is a transposition of a pair of coordinates; this pair is
necessarily contained in only one of the sets $K_k$,
$k=1,\dots,p$. This proves (\ref{produ}).

If $T=\bbT^n$, then $T$ is a product of tori in irreducible
components. Thus, the case $T=\bbT^n$ follows from the above case,
since the assumptions of the theorem hold true for the group $T^0$
if they hold for $T$ in (\ref{inct}) in Example~\ref{maict}.
\end{proof}

\section{Hulls of isotropy orbits of bounded symmetric
domains\label{isohe}}

We start with a preliminary material on hermitian symmetric spaces
following \cite{Wo} but adapting the exposition  to our purpose in
order to be as self contained as possible. For a subset $X$ of a
Lie algebra $\frg$, $\frz(X)=\{z\in\frg:\,[z,X]=0\}$ is the
centralizer of $X$.  Let $G$ be a simple real noncompact Lie group
with a finite center, $K$ be its maximal compact subgroup, and
$\frg,\frk$ be their Lie algebras, respectively. If the center
$\frz=\frz(\frk)$ of $\frk$ is nontrivial, then $\frg$ is called
{\it hermitian}. Then  $\frk=\frz(\frz)$ and $\dim\frz=1$ (note
that $K$ is irreducible in $\frg/\frk$). Let $\frc$ be a Cartan
subalgebra of $\frk$. Then $\frc$ is also a Cartan subalgebra of
$\frg$ and $\frz\subseteq\frc$. There exists $\sfk\in\frz$ such
that $\ad(\sfk)$ has eigenvalues $0,\pm i$ (it is unique up to a
sign; $\ker\ad(\sfk)=\frk$). Then $\ka=e^{\pi\ad(\sfk)}$ is the
Cartan involution which defines the Cartan decomposition
\begin{eqnarray}\label{carta}
\frg=\frk\oplus\frd,
\end{eqnarray}
where $\frk,\frd$ are eigenspaces for $1,-1$, respectively.
Furthermore, $j=\ad(\sfk)$ is a complex structure in $\frd$. This
defines the structure of a hermitian symmetric space of noncompact
type in $D=G/K$. These spaces can be realized as bounded symmetric
domains in $\bbC^n$ with $K=\Aut_0(D)$. Any irreducible bounded
symmetric domain admits such a realization. Let $\De\subseteq
i\frc^*$ be the root system of $\frg^\bbC$. Each $\al\in\De$
corresponds to an $\sll_2$-triple $\sfh_\al,\sfe_\al,\sff_\al$
such that $i\sfh_\al\in\frc$. Thus, $\al(\sfh_\al)=2$,
$[\sfe_\al,\sff_\al]=\sfh_\al$, and
\begin{eqnarray}\label{roots}
[h,\sfe_\al]=\al(h)\sfe_\al,\quad [h,\sff_\al]=-\al(h)\sff_\al
\end{eqnarray}
for all $h\in\frc^\bbC$. We identify $\frc^\bbC$ and
$(\frc^*)^\bbC$ equipping $\frg$ with an $\Ad(K)$-invariant
sesquilinear inner product and normalize it by the condition
\begin{eqnarray}\label{longn}
\max\{|\al|:\,\al\in\De\}=\sqrt2.
\end{eqnarray}
Then short roots must have length $1$ (note that $G_2$ has no real
hermitian form). The set $\De^\vee=\{\sfh_\al:\,\al\in\De\}$ is
the dual root system. The above normalization implies $h_\al=\al$
for long roots and $h_\al=2\al$ for short ones. Since $\ad(h)$,
$h\in\frc$, has eigenvalues $0$ and $\al(h)$, where $\al\in\De$,
we get $\al(i\sfk)=0,\pm1$, i.e., $i\sfk$ is a {\it microweight}
(of $\De^\vee$). For $s=0,\pm1$, set
\begin{eqnarray}\label{defds}
\De_s=\{\al\in\De:\,\al(i\sfk)=s\}.
\end{eqnarray}
Since $\frk\oplus i\frd$ is a compact real form of $\frg^\bbC$ and
$\spann_\bbR\{i\sfh_\al,\sfe_\al-\sff_\al, i(\sfe_\al+\sff_\al)\}$
is the $\su(2)$-subalgebra corresponding to a root $\al\in\De$, we
have
\begin{eqnarray}\label{spand}
\frd=\spann_\bbR\{\sfe_\al+\sff_\al,\
i(\sfe_\al-\sff_\al):\,\al\in\De_1\}.
\end{eqnarray}
Set
$\frs_\al=\spann_\bbR\{i\sfh_\al,\sfe_\al+\sff_\al,i(\sfe_\al-\sff_\al)\}$.
Then $\frs_\al$ is an $\sll(2,\bbR)$-subalgebra of $\frg^\bbC$ and
\begin{eqnarray}\label{singe}
\al\in\De_{\pm1}\quad\Longleftrightarrow\quad\frs_\al\subseteq\frg.
\end{eqnarray}
Let $E$ be a maximal subset of pairwise orthogonal long roots in
$\De_1$. Set
\begin{eqnarray*}
&\sfh=\sum\nolimits_{\al\in E}\sfh_\al,\quad
\sfe=\sum\nolimits_{\al\in E}\sfe_\al,\quad
\sff=\sum\nolimits_{\al\in E}\sff_\al;\\
&\frs=\sum\nolimits_{\al\in E}\oplus\,\frs_\al.
\end{eqnarray*}
Let $\al,\be\in E$, $\al\neq\be$. Since $\al,\be$ are long and
orthogonal, $\pm\al\pm\be\notin\De$.  Hence,
\begin{eqnarray}\label{comms}
\al,\be\in E,\
\al\neq\be\quad\Longrightarrow\quad[\frs_\al,\frs_\be]=0.
\end{eqnarray}
It follows that $\sfh,\sfe,\sff$ is an $\sll_2$-triple and $\frs$
is a subalgebra of $\frg$. Set
\begin{eqnarray}
\thh=e^{\frac14\pi\ad(\sfe-\sff)},\label{defth}\\
\fra=\spann_\bbR\thh E.\label{defa}
\end{eqnarray}
Here is the standard realization of root systems $B_n$ and $C_n$:
\begin{eqnarray*}
&B_n=\{\pm\eps_k\pm\eps_l,\ \pm\eps_m:\,\ k,l,m=1,\dots,n,\
k<l\};\\
&C_n=\{\pm\eps_k\pm\eps_l,\ \pm2\eps_m:\,\ k,l,m=1,\dots,n,\
k<l\}.
\end{eqnarray*}
Then $C_n=B_n^\vee$, but $C_n$ does not satisfy (\ref{longn}).
These systems have microweights; up to the action of the Weyl
group, they are:
\begin{eqnarray*}
B_n:&\eps_1;\\
C_n:&\frac{i}2(\eps_1+\dots+\eps_n).
\end{eqnarray*}
There are no other irreducible root systems which have
microweights and contain roots of different lengths. Also, $B_n$
and $C_n$ have the same Weyl group $BC_n=\bbZ^n_2S_n$.

\begin{lemma}\label{maxab}
The space $\fra$ is a maximal abelian  subspace of $\frd$.
\end{lemma}
\begin{proof}
A straightforward calculation with 2-matrices shows that
$\thh\sfh=\sfe+\sff$. By (\ref{comms}),
\begin{eqnarray}\label{thhef}
\thh\sfh_\al=\sfe_\al+\sff_\al\quad\mbox{\rm for all}\ \al\in E.
\end{eqnarray}
It follows from (\ref{spand}) that $\fra\subseteq\frd$. Moreover,
$\fra$ is abelian due to (\ref{comms}). Set $\Xi=\De\cap E^\bot$.
We claim that
\begin{eqnarray}\label{orcen}
\Xi\subseteq\De_0.
\end{eqnarray}
Indeed, a root in $\De_1\cap\Xi$ must be short. This may happen
only in $B_n$ or $C_n$, since $G_2$ and $F_4$ have no microweights
and other irreducible root systems have no roots of different
lengths. In $B_n$, $\sfk$ is a short root and all other short
roots are orthogonal to $\sfk$. Hence, they do not belong to
$\De_1$. In $C_n$, $E=\{2\eps_1,\dots2\eps_n\}$; then
$\Xi=\emptyset$. Since $\De_{-1}=-\De_1$, this proves
(\ref{orcen}).

Set $\frb=E^\bot\cap\frc$ and
$\frm=\spann_\bbC\{\sfe_\al,\sff_\al:\,\al\in\Xi\}$. It follows
from (\ref{orcen}) that $\frm\subseteq\frk^\bbC$. Clearly,
$\frz(E)=\frc^\bbC\oplus\frm$. The space $\frm$ is
$\thh$-invariant, because $\thh$ fixes roots in $\Xi$. Due to
(\ref{defa}), we get
\begin{eqnarray*}
\frz(\fra)=\thh\frz(E)=\frb^\bbC\oplus\fra^\bbC\oplus\frm
\end{eqnarray*}
Since $\frb^\bbC\oplus\frm\subseteq\frk^\bbC$, this implies
$\frz(\fra)\cap\frd=\fra$.
\end{proof}

The projection of $\thh\De$ into $\fra$ is the {\it restricted
root system} $\De_\fra$ (it is also the set of roots for
$\ad(\fra)$ in $\frg$). The group
\begin{eqnarray*}
W=\{\Ad(g):\,g\in K,\ \Ad(g)\fra=\fra\}|_\fra,
\end{eqnarray*}
acting in $\fra$, is the Weyl group of $\fra$.

In what follows, we denote by $\frv$ the complexification of
$\fra$ with respect to the complex structure $j$ (thus,
$\frv\subset\frd$). The set $\thh E$ is a base in $\frv$;
enumerating it, we identify $\frv$ with $\bbC^n$. Set
$\frt=\spann_\bbR iE$, $T=\exp\frt$, $H=WT$. The torus $T=\bbT^n$
is a maximal compact subgroup in the group $\exp\frs\subseteq G$.

\begin{proposition}\label{pasta}
The following assertions hold:
\begin{enumerate}
\item $\De_\fra$ is a root system of type $BC_n$ or $C_n$; \item
the pair $(\frv,H)$ is standard with $T=\bbT^n$.
\end{enumerate}
\end{proposition}
\begin{proof}
(1). Let $\De_\fra\setminus \thh E$ contain a long root $\al$.
Then $\al=\frac12(\al_1+\al_2+\al_3+\al_4)$ for some
$\al_1,\dots,\al_4\in \thh E$, since $|\al|^2=2$  and
$\scal{\al}{\be}=0,\pm1$ for all $\be\in E$ due to the
normalization (\ref{longn})  (note that $\al,\be$ generate $A_2$
if $\scal{\al}{\be}\neq0$). Roots $\al,\al_1,\dots,\al_4$ generate
$D_4$, since only $A_4$ and $D_4$ among irreducible systems of
rank 4 consist of roots of equal length, but $A_4$ does not
contain an orthogonal base. Since $\scal{i\sfk}{\be}=1$ for all
$\be\in E$, the projection of $i\thh\sfk$ into $\spann_\bbR D_4$
is a microweight $\om$ such that $\scal{\om}{\al_k}=1$,
$k=1,\dots,4$, but $D_4$ has no microweight with this property (in
the realization above, $D_4=B_4\cap C_4$ and the microweights are
either $\pm\eps_k$ or
$\frac12(\pm\eps_1\pm\eps_2\pm\eps_3\pm\eps_4$)). Thus, $E\cup
(-E)$ is the set of all long roots in $\De_\fra$. According to the
classification of irreducible root systems, only $C_n$ and
$BC_n=B_n\cup C_n$ has the property that linearly independent long
roots are mutually orthogonal.

(2). The maximal compact subgroup of the group corresponding to
$\frs$ is $\bbT^n$. Hence, $T=\bbT^n\supset\bbZ_2^n$. Systems
$C_n$ and $BC_n$ have the same Weyl group $W=BC_n$. Therefore,
$H=WT=S_n\bbT^n$.
\end{proof}
Let $D$ be a bounded symmetric domain in a complex linear space
$\frd$ (may be, reducible) and $\frv\subseteq \frd$ be the complex
linear span of a maximal abelian subspace in $\frd$ (thus, we
identify $\frd$ with the corresponding space in the Cartan
decomposition (\ref{carta}), which is induced by the Cartan
involutions in irreducible components). Let $\Aut_{00}(\frv,D)$
denote the subgroup of all linear transformations in $\Aut(D)$
which keep $\frv$ and each irreducible component of $D$.
\begin{corollary}\label{conop}
Let $F$ be a subgroup of $S_n$, $G=F\bbT^n\subset\GL(n,\bbC)$.
Then $G$ satisfies condition (2) of Theorem~\ref{huini} if and
only if $(V,G)$ is isomorphic to a pair
$\left(\frv,\Aut_{00}(\frv,D)\right)$ for a bounded symmetric
domain $D$.
\end{corollary}
\begin{proof}
All pairs $(\bbC^n, S_n\bbT^n)$ appear as
$\left(\frv,\Aut_{00}(\frv,D)\right)$ for matrix balls $D$. It
remains to combine Theorem~\ref{huini} and
Proposition~\ref{pasta}.
\end{proof}
It is possible now to describe hulls of $K$-orbits in $\frd$ (with
respect to the complex structure $j$) it terms of
Proposition~\ref{autdn}. The key point is that $K$ is polar in
$\frd$: each $K$-orbit meets $\fra$ orthogonally (i.e., $\fra$ is
a {\it Cartan subspace}). This is true, since all maximal abelian
subspaces are conjugate in $\frd$ by $K$, $\ad(a)$ is symmetric if
$a\in\frd$ and, for a generic $a\in\fra$, $\ker\ad(a)=\fra$;
hence,
\begin{eqnarray}\label{polar}
[a,\frg]=\fra^\bot.
\end{eqnarray}
We may include the linear base in $\frv$ into a base in $\frd$ as
the first $n$ vectors of the latter. Then $z_1,\dots,z_n$ are
coordinates in $\frv$ and linear functions in $\frd$. The
functions $\mu_k$ in (\ref{mukde}) admit a $K$-invariant extension
to $\frd$:
\begin{eqnarray}\label{extmu}
\mu_k(z)=\sup\{|(gz)_{1}\dots (gz)_{k}|:\,g\in K\},
\end{eqnarray}
where $k=1,\dots,n$. The following lemma shows that (\ref{extmu})
is an extension indeed.
\begin{lemma}\label{coinc}
For $z\in\frv$, {\rm(\ref{mukde})} and {\rm(\ref{extmu})}
coincide.
\end{lemma}
\begin{proof}
It follows from (\ref{polar}) that any critical point of the
linear function $\re z_1$ on the orbit $Kz$ belongs to $\fra$. If
the lemma is not true, then there exist $z\in\frv$ and
$k\in\{1,\dots,n\}$ such that $|(gz)_k|>|z_k|$. Transformations in
$S_n$ and $T$ reduce the problem to the case $z_1>\dots>z_n>0$ and
$k=1$, but then the assumption implies that $\re z_1$ attains its
maximal value on $Kz$ outside of $\fra$.
\end{proof}
\begin{proposition}\label{hulfe}
For any $v\in\frd$, $\wh{Kv}=\{z\in\frd:\,\mu_k(z)\leq\mu_k(v),\
k=1,\dots,n\}$.
\end{proposition}
\begin{proof}
Due to (\ref{extmu}), each $\mu_k$ is a supremum of absolute
values of holomorphic polynomials. Hence, the right-hand side is
polynomially convex. Thus, it includes $\wh{Kv}$. The inverse
inclusion holds, since each $K$-orbit intersects $\frv$ by an
$H$-orbit and hulls of $H$-orbits are distinguished in $\frv$ by
the same inequalities according to Proposition~\ref{autdn} and
Lemma~\ref{coinc}.
\end{proof}
The functions $\mu_k$ can be written in more invariant terms. To
do it, note that the Weyl group of $\De_\fra$ has the form
$\bbZ_2^n S_n$ in the base $\thh E$ by (\ref{defa}); thus,
$z_k=\al_k(z)$, $k=1,\dots,n$, where $\al_k\in\thh E$ and
$z\in\fra$. Therefore, $z_k$ are eigenvalues of $\ad(z)$ in the
subspace generated by the corresponding root vectors. The problem
is to distinguish this subspace (in fact, we use a slightly
different version). After that, functions $\mu_k$ can be defined
as norms of some operators according to the following lemma (this
observation was used in \cite{Ko} in another context).
\begin{lemma}\label{normk}
Let $V$ be a Euclidean space and $A$ be a symmetric nonnegative
operator in $V$ with eigenvalues
$\la_1\geq\la_2\geq\dots\geq\la_m\geq0$, where $m=\dim V$. Let
$A^{\wedge k}$ be its natural extension to the $k$-th exterior
power $V^{\wedge k}=\bigwedge^k V$. Then
\begin{eqnarray*}
\|A^{\wedge k}\|_{V^{\wedge k}}=\la_1\dots\la_k,
\end{eqnarray*}
where $\|\ \|_k$ is the operator norm with respect to the inner
product in $V^{\wedge k}$.
\end{lemma}
\begin{proof}
The norm of a nonnegative symmetric operator is equal to its
maximal eigenvalue.
\end{proof}
Let $v\in\frg$ be semisimple and $\pi(v)$ denote the projection
onto $\ker\ad(v)$ along other eigenspaces of $\ad(v)$ (note that
$\pi(v)$ is a function of $\ad(v)$, since it is the residue at
zero of the resolvent of $\ad(v)$). Set
\begin{eqnarray*}
&a(v)=\ad([v,[v,\sfk]])\pi(v)\ad(\sfk),\label{defop}\\
&p_k(v)=\|a(v)^{\wedge k}\|_{\frg^{\wedge k}},\quad
k=1,\dots,n.\label{defno}
\end{eqnarray*}
The space $\frd$ is $a(v)$-invariant and $\ker a(v)\supseteq\frk$.
We assume that $\frg$ is equipped with some $K$-invariant inner
product, which extends the inner product in $\frd$. It follows
from the calculation below that $a(v)$ is symmetric and has range
$\ad(\sfk)\fra$. Let $n=\dim\fra$ be the rank of the symmetric
space $D$. It is equal to the codimension of a generic $K$-orbit
in $\frd$.
\begin{theorem}\label{last}
For any $v\in\frd$,
\begin{eqnarray*}
\wh{Kv}=\{z\in\frd:\,p_k(z)\leq p_k(v),\ k=1,\dots,n\}.
\end{eqnarray*}
\end{theorem}
\begin{proof}
It is sufficient to prove the assertion for a generic $v\in\frd$.
Clearly, $p_k$ are $K$-invariant. Hence, we may assume $v\in\fra$.
Then, by (\ref{defa}) and (\ref{thhef}),
\begin{eqnarray*}
v=\sum\nolimits_{\al\in E} v_\al(\sfe_\al+\sff_\al),
\end{eqnarray*}
where $ v_\al\in\bbR$. According to (\ref{roots}) and
(\ref{defds}), $[\sfk,v]=\sum\nolimits_{\al\in E}
iv_\al(\sfe_\al-\sff_\al)$. Thus,
\begin{eqnarray*}
[v,[v,\sfk]]=\sum\nolimits_{\al\in E} 2iv_\al^2\sfh_\al
\end{eqnarray*}
due to (\ref{comms}). Also, (\ref{comms}) implies that
$\ad(i\sfh_\al)$ keeps $\frv$ and has eigenvalues $0,\pm2i$ in it
for each $\al\in E$. Therefore, $\frv$ is
$\ad([v,[v,\sfk]])$-invariant and its eigenvalues are
$\pm4v_\al^2i$, $\al\in E$. Since $\ad(\sfk)\frg=\frd$ and
$\pi(v)\frd=\fra$ for a generic $v\in\fra$, the space $\frv$ is
$a(v)$-invariant; moreover, $a(v)\frg=\ad(\sfk)\fra\subseteq\frv$.
Thus, $a(v)$ has eigenvalues $0,\pm4v_\al^2$ in $\frg$. According
to Lemma~\ref{normk} and (\ref{mukde}),
\begin{eqnarray}\label{pkmuk}
p_k(v)=4^k\mu_k^2(v)
\end{eqnarray}
for $v\in\frv$ and $k=1,\dots,n$. Since $p_k$ and $\mu_k$ are
$K$-invariant, (\ref{pkmuk}) holds for all $v\in\frd$. The theorem
follows from Proposition~\ref{hulfe}.
\end{proof}
\begin{corollary}
Functions $p_k$, $k=1,\dots,n$, are plurisubharmonic in $\frd$
with respect to the complex structure $j=\ad(\sfk)$.
\end{corollary}
\begin{proof}
By (\ref{pkmuk}) and (\ref{extmu}),
\begin{eqnarray*}
p_k(z)=4^k\sup\{|(gz)_{1}^2\dots (gz)_{k}^2|:\,g\in K\}.
\end{eqnarray*}
The right-hand side is plurisubharmonic, since the functions
$z_k^2$ are $j$-holomorphic and $j$ is $K$-invariant.
\end{proof}
One can get the same functions $p_k$ by replacing $\frg$ with
$\frd$, endowed with the complex structure $j$, and $a(v)$ with
$\ad([v,jv])(\pi(v)+\pi(jv))$.

\subsection*{Acknowledgement} I thank Anton Pankratiev for helpful comments.

\bibliographystyle{amsalpha}

\begin{thebibliography}{A}
\bibitem{An}
J.T. Anderson, \textit{Polnomial hulls of sets invariant under an
action of the special unitary group}, Can. J. Math., v.
\setcounter{xxx}{40}{\Roman{xxx}}, (1988) no. 5, 1256-1271.
\bibitem{Bj}
J.-E. Bj\"ork, \textit{Compact groups operating on Banach
algebras}, Math. Ann., v. 205 (1973), no. 4, 281--297.
\bibitem{Da}
J. Dadok, \textit{Polar coordinates induced by actions of compact
Lie groups}, Transactions Amer. Math. Soc., 288  (1985), 125-137.
\bibitem{DG}
A. Debiard, B. Gaveau, \textit{Equations de Cauchy–Riemann sur
SU(2) et leurs enveloppes d'holomorphie},  Can. J. Math., v. 38
(1986), 1009-1024.
\bibitem{Ga}
T. Gamelin, \textit{Uniform Algebras}, Prentice-Hall, Englewood
Cliffs, N. J., 1969.
\bibitem{Gi}
V.M. Gichev, \textit{Maximal ideal spaces of invariant function
algebras on compact groups}, preprint, 46 pp., available at
http://arXiv.org/abs/math/0603449
\bibitem{GL}
V.M. Gichev, I.A. Latypov, \textit{Polynomially Convex Orbits of
Compact Lie Groups}, Transformation Groups, v. 6 (2001) no. 4,
321-331.
\bibitem{FB}
J. Faraut, L. Bouattour, \textit{Enveloppes polynˆomiales
d'ensembles compacts invariants}, Math. Nachrichten 266 (2004),
20-26.
\bibitem{Ka}
W. Kaup,\textit{ Bounded symmetric domains and polynomial
convexity}, Manuscr. Math. 114 (2004), No.3, 391-398.
\bibitem{KZ}
W. Kaup, D. Zaitsev, \textit{ On the CR-structure of compact group
orbits associated with bounded symmetric domains}, Invent. Math.,
153 (2003), 45-104.
\bibitem{Kan}
J. Kane, \textit{Maximal ideal spaces of U-algebras}, Illinois J.
Math., v.27 (1983), No 1, 1-13.
\bibitem{Ko}
B. Kostant, \textit{On convexity, the Weyl group and the Iwasawa
decimposition}, Ann. Sci. Ecole. Norm. Super., 6 (1973), 413--455.
\bibitem{VP}
V.L.~Popov,  E.B.~Vinberg, \textit{Invariant theory}, Itogi Nauki
i Tekhniki, Sovr. Probl. Mat. Fund. Napravl., v.~55, VINITI,
Moscow 1989, pp.~137--309 (in Russian); English transl.: Algebraic
Geometry IV, Encyclopaedia of Math. Sciences, v.~55,
Springer-Verlag, Berlin 1994, pp.~123-278.
\bibitem{Sa}
C. Sacr´e, \textit{Enveloppes polynomiales de compacts}, Bull.
Sci. math. 116 (1992), 129-144.
\bibitem{Wo}
J.A. Wolf, \textit{Fine Structure of Hermitian Symmetric Spaces},
Symmetric spaces (Short Courses, Washington Univ., St. Louis, Mo.,
1969-1970), Pure Appl. Math., Vol. 8. New York: Dekker 1972,  pp.
271-357.
\end{thebibliography}

\end{document}